

\baselineskip=14pt
\parskip=10pt

\font\eightrm=cmr8 

\magnification=\magstephalf
\def\N{{\cal N}}
\def\P{{\cal P}}
\def\Q{{\cal Q}}
\def\R{{\cal R}}
\def\1{{\overline{1}}}
\def\2{{\overline{2}}}
\def\Tilde{\char126\relax}
\parindent=0pt
\overfullrule=0in
\def\frac#1#2{{#1 \over #2}}
\bf
\centerline
{
The Quasi-Holonomic Ansatz and Restricted Lattice Walks
}
\rm
\bigskip
\centerline 
{ 
{\it Manuel KAUERS}\footnote{$^1$}
{
{\eightrm
Supported in part by the Austrian Science Foundation (FWF),
grants SFB F1305 and P19462-N18. 
}
}
and {\it Doron ZEILBERGER}\footnote{$^2$}
{
{\eightrm
Supported in part by the USA National Science Foundation.
}
}
}

\bigskip
\qquad\qquad\qquad\qquad
{\it Dedicated to Gerry Ladas on his 70th Birthday}
\bigskip

{\bf Preface: A One-Line Proof of Kreweras' Quarter-Plane Walk Theorem}

See:  {\tt http://www.math.rutgers.edu/{\Tilde}zeilberg/tokhniot/oKreweras} .

{\bf Comments}: The great enumerator
Germain Kreweras empirically  discovered this
intriguing fact, and then needed lots of pages[K], and
lots of human ingenuity, to prove it. Other great enumerators,
for example, Heinrich Niederhausen[N], Ira Gessel[G1],
and Mireille Bousquet-M\'elou[B], found other ingenious,
``simpler'' proofs. Yet none of them is as simple as ours!
Our proof (with the generous help of our faithful computers)
is ``ugly'' in the traditional sense, since it would be
painful for a lowly human to follow all the steps.
But according to {\it our} humble aesthetic taste, this proof is much
more elegant, since it is (conceptually) {\it one-line}. So what if
that line is rather long (a huge partial-recurrence equation
satisfied by the general counting function), it occupies less
storage than a very low-resolution photograph.

{\bf Unrestricted Lattice Walks}

Suppose that you are walking, in the $d$-dimensional hyper-cubic
lattice $Z^d$, starting at the origin, and
at each time-unit (you can call it a nano-second if you are a fast-walker,
or a year if you are slow), you are allowed to use {\it any} step
from a certain {\bf finite} set of {\bf fundamental steps}
$$
S=\{ (s_1, \dots, s_d) \} \quad,
$$
where each fundamental step can have {\it arbitrary} integer components
(i.e. negative, positive, or zero).

For example, for the simple lattice (``random'') 
walk on the line, we have $S=\{-1, 1\}$,
while for the simple random walk on the
two-dimensional square lattice, 
we have $S=\{(1,0),(-1,0),(0,1),(0,-1)\}$.
To cite another example, a
Knight, on an infinite chessboard, is allowed any of the
following eight steps:
$$
S=\{ \, (\pm 2, \pm 1) \, ,  \, (\pm 1, \pm 2) \,  \} \quad .
$$

The quantity of interest is the $d+1$-variable discrete function,
let's call it
$$
F(m;n_1, \dots, n_d) \quad,
$$
that counts the number of ways of walking from the origin $(0, \dots, 0)$
to the point $(n_1, \dots, n_d)$ in {\it exactly} $m$ steps.

Often, one is interested, more specifically, in
$f(m):=F(m; 0, \dots, 0)$, the number of such walks that
return to the origin after $m$ steps, and, of course we have
that $g(m)$, the {\it total} number of walks with $m$ steps, 
at the present {\it unrestricted} case, is trivially
$\vert S \vert^m$.

It is very easy to write down the (full) {\it  generating function}
of $F$:
$$
\tilde F (t;x_1, \dots, x_d) \, := \,
\sum_{m=0}^{\infty} \sum_{n_1=-\infty}^{\infty} \dots
\sum_{n_d=-\infty} ^{\infty}
F(m;n_1, \dots, n_d) \, t^m {x_1}^{n_1} \dots {x_d}^{n_d}  \quad .
$$
Indeed, the readers will have no trouble convincing themselves that
$$
\tilde F (t;x_1, \dots, x_d)=
\left [
1-t \left ( \sum_{(s_1, \dots , s_d) \in S} 
{x_1}^{s_1} \dots {x_d}^{s_d} \right )
\right ]^{-1} \quad,
$$
which is a {\it rational} function of its variables, and
it should be interpreted as a formal power series in
$t$ whose coefficients are {\it Laurent polynomials} in
$(x_1, \dots, x_d)$.

It follows immediately from ``general holonomic nonsense'' 
[Z1][WZ] that $F(m; n_1, \dots, n_d)$ is {\it completely} holonomic,
i.e. it satisfies $d+1$ {\it pure} (homogeneous) linear recurrences
with polynomial coefficients, one for each of its
arguments. (Generically speaking. In some degenerate cases
some or all of these $d+1$ equations coincide, and one needs
more equations to describe the function.).

More verbosely,
there exists a positive integer L, and polynomials 
$$
p_0(m;n_1, ..., n_d) \, , \, p_1(m;n_1, ..., n_d) \, , \, \dots 
\, , \, p_L(m;n_1, ..., n_d)  \quad  ,
$$
such that
$$
\sum_{i=0}^{L} p_i(m; n_1, \dots, n_d) F(m+i \, ; \, n_1, \dots, n_d) 
\quad =0 \quad ,
$$
for {\it all} $m \geq 0$ and $(n_1, \dots, n_d) \in Z^d$,
and for each dimension $n_i$, ($i=1 \dots d$), there exists
a positive integer $K_i$, and polynomials 
$q^{(i)}_j (m; n_1, \dots, n_d)$, $j=0 \dots K_i$, such that
$$
\sum_{j=0}^{K_i}  
q^{(i)}_j (m; n_1, \dots, n_d)
F(m \, ; \, n_1, n_{i-1}, n_i+j, n_{i+1} \dots, n_d) 
\, = \, 0 \quad .
$$

Furthermore, thanks to [MZ]
(that contains, among other things, a  multi-variable
extension of the Almkvist-Zeilberger[AZ] algorithm), and [W], one can
actually {\it explicitly} find these recurrences. 
However, for $d$ larger
than $4$ and/or for large sets $S$, it soon becomes impractical
with today's computers.

\vfill\eject

{\bf Restricted Lattice Walk}

Very often, in real life, we would like to stay in certain
sub-regions of $Z^d$. In this case, it is no longer true that
the counting function $F$ is necessarily holonomic, as shown by
Mireille Bousquet-M\'elou and Marko Petkovsek in a seminal paper
[MP]. But, sometimes it is still holonomic, because of the
``nice'' structure of the restricted region.

For example, if the the set of steps, $S$, consists of
the unit positive steps in $d$ dimensions, and one is only allowed
to stay in $n_1 \geq n_2 \geq \dots \geq n_d $, we have
the famous $d$-dimensional {\it ballot} problems, that
is equivalent to the problem of enumerating standard Young Tableaux of
shape $(n_1, \dots , n_d)$. Here we famously have
the Young-Frobenius-MacMahom formula, that
$f(n_1, \dots, n_d):=F(n_1 + \dots + n_d; n_1, \dots, n_d)$
is given by
$$
f(n_1, \dots, n_d)
\, = \,
\prod_{1 \leq i < j \leq d} (n_i-n_j+j-i) 
\cdot 
{{(n_1 + \dots + n_d)!}
\over
{(n_1+d-1)! (n_2+d-2)! \dots  (n_d)!}} \quad,
$$
that immediately implies that not only is it holonomic, but
the relevant recurrences for $f$ are {\it first-order}
in each of its variables, since $f$ is
expressible in {\it closed-form}.

There are other examples, even allowing negative steps,
where one still stays in the holonomic realm,
see for example [GZ]. This happens because 
if you put mirrors on the bounding hyper-planes, the 
group generated by the reflections is finite (the so-called Weyl, or Coxeter group),
and the set of steps is invariant under that group.

{\bf Kreweras' Walks}

But things start to get complicated very soon. Consider
the following set of {\bf three} steps
$$
S=\{ (-1,0),(0,-1),(1,1) \} \quad ,
$$
walking in {\bf two} dimensions, and staying in the
{\bf positive quadrant}, i.e. one must stay in 
the region $\{ (n_1,n_2) \, | \, n_1 \geq 0 \, , \, n_2 \geq 0 \}$.

Obviously, $F(m; n_1, n_2)$, 
defined for $m \geq 0, n_1 \geq -1,n_2 \geq -1$,
satisfies the following
simple recurrence
$$
F(m; n_1,n_2)\, =\, 
F(m-1; n_1+1,n_2)+F(m-1; n_1,n_2+1)+F(m-1; n_1-1,n_2-1) \quad,
$$
(whenever $m \geq 1$ and $n_1,n_2 \geq 0$),

subject to the {\bf initial condition}:
$$
F(0; n_1, n_2)= 
\cases{ 1
  ,& if ($n_1$,$n_2$)=(0,0) , \cr
0, & otherwise .\cr}
$$
and the {\bf boundary conditions}:
$$
F(m; n_1,n_2) \, = 0 \, \, if \,\, n_1=-1 \, \, \, or \,\,\, n_2 =-1 \quad .
$$

Surprisingly, $F(m; 0,0)$ is {\bf closed-form}. In a classic
paper, Germain Kreweras[K] proved that 
$$
F(3n; 0,0) \, = \, { {4^n} \over {(n+1)(2n+1)} }
{{3n} \choose {n}} \quad ,
$$
(of course $F(m;0,0)=0$ if $m$ is not a multiple of $3$).

A naive approach would be to try and conjecture a closed-form
formula, in terms of $m,n_1,n_2$,
for the general $F(m; n_1,n_2)$, verify that this formula
obeys the above simple recurrence and the initial and boundary
conditions, and finally plug-in $n_1=0,n_2=0$.

Alas, while $F(m;0,0)$ is almost as nice as could be,
the general $F(m;n_1,n_2)$ is a huge mess, and the above
approach is doomed to failure, at least if taken literally.
We will later show how to rescue this simple-minded approach,
by reasoning in the holonomic
(or if necessary, {\it quasi-holonomic}) realm.

{\bf Approaches}

The most successful approach so far, was to derive
a {\bf functional equation} for the generating function, 
using combinatorial ([K]) or
probabilistic ([G1]) reasoning, followed sometimes by the
Kernel method, brought to new heights by ``La Mireille''[B].
A very nice systematic study of the successes of the Kernel method,
still in the quarter plane, and with {\bf exactly} {\it three}
steps, all with coordinates between $-1$ and $1$, was undertaken
by Marni Mishna[Mi].

{\bf Ira Gessel's Intriguing Conjecture}

If the set of steps is
$$
S=\{ (-1,0),(1,0),(-1,-1), (1,1) \} \quad ,
$$
still staying in the positive quadrant ($\{(x,y) \, \vert \, x \geq 0, y\geq 0 \}$),
then Ira Gessel[G2] discovered empirically that
(recall that $(a)_n:=a(a+1) \cdots (a+n-1)$),
$$
F(2n; 0,0)= 16^{n} { {(5/6)_n (1/2)_n} \over {(2)_n (5/3)_n}} \quad .
$$
(Of course F(2n+1;0) =0).
At this time of writing, as far as we know, this remains unproved.
The Kernel method, so far, did not succeed,
perhaps because that {\it now} there are {\bf four} steps.

{\bf The Holonomic Approach}

In [Z2], an empirical-yet-rigorous approach for enumerating
{\it unrestricted} lattice paths was suggested, using the
{\it holonomic} ansatz. 
This method should, (and indeed does!, see the output 
at this article's webpage)
succeed in doing the Kreweras problem. 
We now know, a posteriori, that the
full generating function $\tilde F(t;x_1,x_2)$ for
Kreweras walk is even {\it algebraic}, and hence
{\it a fortiori}, holonomic. Hence, there exists, a (giant!)
linear recurrence operator
$$
\P(M,m,n_1,n_2) \quad,
$$
where $M$ is the shift operator in $m$,
(i.e. $Mf(m):=f(m+1)$ for any function $f(m)$ ),
annihilating $F(m;n_1,n_2)$.
It turns out that $\P(M,m,n_1,n_2)$ is extremely complicated,
but once found,  plugging-in $n_1=0,n_2=0$ gives an operator,
$\P(M,m,0,0)$, annihilating $F(m;0,0)$.

\vfill\eject

{\bf How to prove that the empirically-derived operator does
indeed annihilate $F$?}

Let's restrict attention to the quarter-plane. Similar
reasonings apply to higher dimensions and more general
regions.

Given a set of steps $S$, our discrete function
$F(m; n_1,n_2)$ satisfies the recurrence
$$
F(m;n_1,n_2)= \sum_{(s_1,s_2) \in S} F(m-1;n_1-s_1,n_2-s_2) \quad,
$$
which means that $F(m; n_1, n_2)$ is {\bf annihilated} by the
linear recurrence operator with {\bf constant} coefficients
$$
\Q=
1-M^{-1} \left ( \sum_{(s_1,s_2) \in S} N_1^{-s_1}N_2^{-s_2} \right ) \quad .
$$

We want to prove that $\Q F=0$ plus the obvious
initial and boundary conditions, imply that $\P F=0$.

Let's call an operator {\it good} it it only contains non-negative exponents
of the shift operators $N_1,N_2$. For example $1-M^{-1}N_1^2-M^{-2}N_2$ is good
but $1-M^{-1}N_1^{-1}$ is not.

By taking commutators, or otherwise, we find,
calling $\P_0=\P$, a
{\bf sequence} of {\bf good operators}, 
$\R_0(m,n_1,n_2,M,N_1,N_2),  \dots , \R_d(m,n_1,n_2, M,N_1,N_2)$, and operators  
$\P_1(m,n_1,n_2,M,N_1,N_2)$,  $\dots$ , $\P_d(m,n_1,n_2,M,N_1,N_2)$,
of lower-and-lower degrees such that
$$
\Q \P_0= \R_0 \Q + \P_1 \quad ,
$$
$$
\Q \P_1= \R_1 \Q + \P_2 \quad ,
$$
$$
\dots
$$
$$
\Q \P_d= \R_d \Q + \P_{d+1} \quad ,
$$
with $\P_{d+1}=0$. Since $\R_d$ is ``good'',
and since $\Q F =0$, we have that $R_d \Q F=0$ and
hence $\Q [\P_d F]=0$. Then check that the boundary conditions
for $\P_d F$ are the same and the initial condition is
identically $0$ to deduce that $\P_d F=0$. By backwards induction
we (or rather our computer, it can all be mechanized) in turn, proves
$\P_{d-1}F=0 \, , \, \P_{d-2}F=0 \, , \, \dots \, , \,  \P_{0} F=0$.

Note that if you don't insist that the $\R_i$'s are ``good'' one can
always take $\R_i=\P_i$, and $\P_{i+1}$ is simply the commutator
of $\Q$ and $\P_i$, for $i=0, 1, \dots, d$. Since $\Q$ is 
constant-coefficients, taking commutators with
any operator with polynomial coefficients, always decreases the
degree of the polynomial coefficients, so if the degree is $d$,
eventually, after $d+1$ iterations,
we get that $\P_{d+1}=0$. If we want the $\R_i$
to be good, we have to adjust things to be good.

In fact, for the lattice-paths-counting problems  treated here,
with the time
variable $m$, starting at time $m=0$ at the origin, 
it is not really necessary to demand that the $\R_i$ be
``good''. We can consider the function $F$ to be defined everywhere,
with $0$ at the forbidden region, and rephrase that
$\Q F=\delta(n_1,n_2)$ when $m=0$, where $\delta(n_1,n_2)$ is
the discrete delta function that is $1$ at the origin and
$0$ elsewhere.

{\bf The Quasi-Holonomic Approach} 

For the sake of exposition, let's stay in the plane
(analogous reasoning applies in general),
and let $n$ denote discrete time and $(a,b)$ discrete space.

As mentioned above,
Mireille Bousquet-M\'elou and Marko Petkovsek
proved that it is not always true, for arbitrary steps
and arbitrary boundaries, that the counting function is
holonomic. It is probably usually false, and 
the holonomicity of the
Kreweras walks, and the few other cases in which it may hold,
are just flukes (or follow from other considerations).

But who cares about holonomicity? Maybe it is asking way too much.
Suppose, like, in the case of Gessel's conjecture mentioned above,
$F(n;0,0)$ turns out to be holonomic in the {\it single} variable $n$.

If $F(n;a,b)$ is holonomic in {\it all} its arguments, then
there exist {\it three} independent, pure recurrence operators
$$
\P_1(n,a,b, N) \quad,  \quad \P_2(n,a,b, A) \quad,  \quad
\P_3(n,a,b, B) \quad ,
$$
annihilating $F$. In particular, $\P_1(n,0,0,N)$
would give us the desired operator.

But, very likely, $F(n;a,b)$ is {\it not}
holonomic, and even if it is, like in Kreweras' case,
$\P_1(n,a,b,N)$ is too big. What do we do now?
Something much more modest would do the job!

All we need is {\it one} linear recurrence operator with
polynomial coefficients of the form
$$
\R(a,b,n,A,B,N)=
\R_0(n,N)+ a \R_1(a,b,n,A,B,N)+b \R_2(b,n,A,B,N) \quad ,
$$
with  $\R_0 \neq 0$.

Once found, empirically, one can prove that it 
annihilates our counting function $F(n,a,b)$ as above,
by constructing a sequence of operators (by taking
commutators, and possibly tweaking to get good operators).
Once $\R$ is found, and proved to indeed annihilate
$F(n,a,b)$ (all of which should be done completely
automatically by the computer), all we have to do is
plug-in $a=0, b=0$ in 
$$
\R(a,b,n,A,B,N) F(n,a,b)=0 \quad ,
$$
and get that
$$
\R_0(n,N) F(n,0,0)=0 \quad.
$$
QED.

Our ``one-line'' proof of Kreweras' theorem, mentioned in the
preface, used this quasi-holonomic ansatz, even though, in this case,
it is known that the counting function is holonomic. Staying
within the holonomic ansatz would have made the ``one-line'' yet
longer 
and its computation yet slower.
(A holonomic operator $\R(a,b,n,N)$ for the Kreweras walks is, 
for comparison, also available from the webpage of this article:
{\tt krewerasComplete}).

Analogs to Kreweras' theorem can be found effortlessly for all the 
eleven walks that Mishna[Mi] has isolated as being essentially different. 
The results are as follows:

\medskip
\halign{\indent\hfil#&\quad#\hfil&\quad#\hfil\cr
   & step set & number of closed paths \cr
 1 & $\{(0,1),(1,1),(1,0)\}$ & $f(n,0,0)=0$ \cr
 2 & $\{(0,1),(1,1),(-1,-1)\}$ & $f(2n,0,0)=\frac{4^n(1/2)_n}{(1)_{n+1}}$ \cr
 3 & $\{(0,1),(1,1),(1,-1)\}$ & $f(n,0,0)=0$ \cr
 4 & $\{(0,1),(0,-1),(1,-1)\}$ & $f(2n,0,0)=\frac{4^n(1/2)_n}{(1)_{n+1}}$ \cr
 5 & $\{(-1,0),(0,-1),(1,1)\}$ & $f(3n,0,0)=\frac{2\cdot27^{n-1}(4/3)_{n-1}(5/3)_{n-1}}{(5/2)_{n-1}(3)_{n-1}}$ \cr
 6 & $\{(0,1),(1,0),(-1,-1)\}$ & $f(3n,0,0)=\frac{2\cdot27^{n-1}(4/3)_{n-1}(5/3)_{n-1}}{(5/2)_{n-1}(3)_{n-1}}$ \cr
 7 & $\{(-1,0),(0,1),(1,-1)\}$ & $f(3n,0,0)=\frac{27^{n-1}(4/3)_{n-1}(5/3)_{n-1}}{(3)_{n-1}(4)_{n-1}}$ \cr
 8 & $\{(0,1),(-1,-1),(1,-1)\}$ & $f(4n,0,0)=\frac{2\cdot64^{n-1} (5/4)_{n-1} (3/2)_{n-1}(7/4)_{n-1}}{(2)_{n-1}(5/2)_{n-1}(3)_{n-1}}$ \cr
 9 & $\{(-1,0),(1,1),(1,-1)\}$ & $f(4n,0,0)=\frac{2\cdot64^{n-1} (5/4)_{n-1} (3/2)_{n-1}(7/4)_{n-1}}{(2)_{n-1}(5/2)_{n-1}(3)_{n-1}}$ \cr
 10 & $\{(-1,1),(0,1),(1,-1)\}$ & $f(n,0,0)=0$ \cr
 11 & $\{(-1,1),(1,1),(1,-1)\}$ & $f(n,0,0)=0$ \cr
}\medskip 

Computer-generated proofs for the non-zero entries can be found
at the webpage of this article.

We have also searched for an operator $\R(a,b,n,A,B,N)$ that would yield a proof
of Gessel's conjecture, but it has turned out that no such operator can be found
whose degree in $A,B,N$ individually is at most~$8$ and whose total degree in
$a,b,n$ is at most~$6$.

{\bf A more refined counting}

Another interesting problem is as follows. Given a  set 
of steps $S=\{ S_i \vert i=1 \dots r \}$,
count the number of
walks with exactly $A_i$ steps of kind $S_i$.
Now the condition that it stays in the quarter-plane
(or half-line, or eighth-space, or whatever), can be expressed
as walks, with {\it positive unit steps} in $\N^r$
confined to the positive sides of certain hyperplane.
For example, for Kreweras's walks, if
$f(a,b,c)$ is the number of walks using $a$ steps of
kind $(-1,-1)$, $b$ steps of kind $(1,0)$ and $c$ steps
of kind $(0,1)$, we are counting walks from the origin
to $(a,b,c)$ staying in $c \geq a$ and $b \geq a$.
Then $f(n,n,n)$ is what we called above $F(3n;0,0)$.

To get the quantity of interest in Gessel's conjecture,
we need to compute
$$
G(n):= \sum_{a=0}^{n} f(a,a,n-a,n-a) \quad .
$$
Even though $f(a,b,c,d)$ is
unlikely to be holonomic, let's hope that it is
quasi-holonomic enough to guarantee that $G(n)$ 
is holonomic in the single variable $n$,
a fact that we already know empirically, but it would be nice to prove it.

The Maple package {\tt WalkCarefully} counts walks this way.

{\bf Open problem (even empirically) }

Is the analog of Kreweras' walk in three dimensions  holonomic?

In other
words  does the sequence $a(n):=$
the number of ways of walking in the positive eigth-space
( $\{(x,y,z) \, \vert \, x \geq 0, y \geq 0, z \geq 0 \}$),
starting at the origin, walking $4n$ steps, and 
returning to the origin, only employing the
steps
$$
\{ \, (-1,-1,-1)\, , \, (1,0,0) \, , \,  (0,1,0) \, , \, (0,0,1) \, \} \quad,
$$
a solution of a linear recurrence equation with polynomial coefficients?
According to our computations, it would satisfy a recurrence
of very high order and degree if it were.

{\bf What About Gessel's Problem}

The {\it holy grail} for lattice-walk-counters, currently,
is a proof of Ira Gessel's conjecture. We strongly believe that
the counting function is {\it quasi-holonomic}, so the present
approach should, {\it at least in principle}, prove it.
But, of course, it remains to be seen whether our proverbial
margin is wide enough to contain the proof.

We also strongly believe that
there is a much simpler proof (in all senses of the word) of that conjecture,
that requires less that 1K of memory. That simple proof would come once
the {\it right} and {\it natural} ansatz to which the 
(restricted) counting function belongs to,
will be discovered. To give an analogy, we can routinely prove that
$$
\sum_{k=0}^{10000000} {{10000000} \choose {k}} x^k y^{10000000-k} \, = \,
(x+y)^{10000000} \quad,
$$
by staying in the {\it polynomial} ansatz. But it would be
much more efficient to first prove that
$$
\sum_{k=0}^{n} {{n} \choose {k}} x^k y^{n-k} \, = \,
(x+y)^{n} \quad,
$$
for {\it all} $n$,
by working in the {\it holonomic} ansatz, using WZ theory, say,
and then, simply, plug-in $n=10000000$.

The hard part, of course, for which we still need humans, is
to {\it chercher l'ansatz}.

{\bf Maple and Mathematica Packages}

This article is
accompanied by four very basic Maple packages, that compute
the counting functions and empirically guess recurrences.

These are:

{\tt HalfLine , OneDimWalks, QuarterPlane, WalkCarefully } .

There are also Mathematica packages 

{\tt Guess , Walks }

All these are available at the webpage of this article:

{\tt http://www.math.rutgers.edu/\Tilde zeilberg/mamarim/mamarimhtml/quasiholo.html}  .

It is hoped that these can be extended to prove, fully automatically,
Gessel's conjecture, as well as make up their own conjectures and proofs
for other sets of steps.

{\bf References}

[AZ] G. Almkvist and D. Zeilberger,
{\it The method of differentiating under the
integral sign}, J. Symbolic Computation {\bf 10}, 571-591 (1990).

[B] M. Bousquet-M\'elou, {\it Walks in the quarter plane:
Kreweras' algebraic model}, Annals of Applied Probability
{\bf 15} (2005), 1451-1491.

[BP] M. Bousquet-M\'elou and M. Petkovsek,
{\it Walks confined in a quadrant are not always D-finite},
Theor. Computer Sci. {\bf 307} (2003), 257-276.

[G1] I. Gessel, {\it A probabilistic method for lattice
path enumeration}, J. Stat. Plann. Inference {\bf 14} (1986),
49-58.

[G2] I. Gessel , {\it private communication}.

[GZ] I. Gessel and D. Zeilberger, {\it Random Walk in a Weyl chamber}, Proc.
Amer. Math. Soc. {\bf 115}, 27-31 (1992).

[K] G. Kreweras, {\it Sur une class de probl\`emes li\'es au
treillis des partitions d' entiers},
Cahiers du B.U.R.O, {\bf 6} (1965), 5-105.

[Mi] M. Mishna, {\it Classifying lattice walks restricted
to the quarter plane}, preprint. Available from
{\tt arxiv.org} (\# 0611651).

[MZ] M. Mohammed and D. Zeilberger,
{\it Multi-Variable Zeilberger and Almkvist-Zeilberger Algorithms and the
Sharpening of Wilf-Zeilberger Theory}, Adv. Appl. Math. {\bf 37}
(2006), 139-152.

[N] H. Niederhausen, {\it The ballot problem with three
candidates}, European J. Combin., {\bf 1} (1980), 175-188.

[W] K. Wegschaider. {\it ``Computer Generated Proofs of Binomial Multi-Sum
Identities}, RISC, J. Kepler University. Diploma Thesis. May 1997.

[WZ] H.S. Wilf and D. Zeilberger,
{\it An algorithmic proof theory for hypergeometric
(ordinary and "q") multisum/integral identities}, Invent. Math. 
{\bf 108}, 575-633 (1992).

[Z1] D. Zeilberger,
{\it A Holonomic systems approach to special functions
identities}, J. of Computational and Applied Math. {\bf 32},
321-368 (1990).

[Z2] {\it The Holonomic ansatz I. Foundations and Applications
to Lattice Paths Counting}, Annals of Combinatorics
{\bf 11} (2007), 227-239.

\bigskip

{\bf Manuel Kauers}, Research Institute for Symbolic Computation,
J. Kepler Univ. Linz, Austria. mkauers@risc.uni-linz.ac.at .

\medskip

{\bf Doron Zeilberger}, Mathematics Department, Rutgers University
(New Brunswick), Piscataway, NJ, USA.  zeilberg@math.rutgers.edu .

\bigskip
First Written: Dec. 5, 2007

This version:  Dec. 14, 2007
\end